\documentstyle[12pt, amssymb, amsmath]{article}
\begin{document}
\newcommand{\bref}[1]{$\mbox{(\ref{#1})}$}
\newcommand{\p}{{\Bbb {R^+}} \rightarrow {\Bbb {C}}}
\newcommand{\r}{{\Bbb {R}} ^d}
\newcommand{\q}{{\mathcal Q}}
\newcommand{\z}{{\Bbb {Z}} ^d}
\newcommand{\s}{\Sigma}
\newtheorem{thm}{Theorem}
\newtheorem{lem}{Lemma}
\newtheorem{conj}{Conjecture}
\newtheorem{rem}{Remark}
\newtheorem{cor}{Corollary}

\title
{{\bf The Uncertainty Principle for certain densities
}}

\author{Oleg Kovrizhkin \thanks{This research was partially partially supported by NSF grant 
DMS 0201099}}

\date{}

\maketitle

\begin{center}
\rm
MIT, 2-273\\
Math Dept\\
77 Mass Ave\\
Cambridge, MA 02139, USA\\
{\sl E-mail address}: oleg@math.mit.edu\\
\end{center}

\begin{abstract} We prove a new version of the Uncertainty Principle of the form $\int |f|^2 
\lesssim \int_{E^c} |f|^2 + \int_{\Sigma ^c}|\hat f|^2 $ where 
the sets $E$ and $\Sigma$ are $\epsilon$-thin in the following sense: $|E \cap D(x, \rho_1(x))| 
\le \epsilon |D(x, \rho_1(x))|$ and $|\Sigma \cap D(x, 
\rho_2(x))| \le \epsilon |D(x, \rho_2(x))|$. This is an intermediate result between 
Logvinenko-Sereda's and Wolff's versions of the Uncertainty Principle.
\end{abstract}

\section{Introduction}\label{intro}

The Uncertainty Principle in Fourier analysis is a statement that a function and its
Fourier transform can not both be concentrated on small sets. Many examples of this
principle can be found in the book by Havin and Joricke \cite{HJ} and in the paper by
Folland and Sitaram \cite{FS}. We will be interested in the following type of the Uncertainty 
Principle for $f \in L^2(\r)$:
\begin{eqnarray}
\int |f|^2 \le C (\int_{E^c} |f|^2 + \int_{\Sigma^c} |\hat f|^2)\label{UPIN}
\end{eqnarray}
where 
$E$ and $\Sigma$ are  \lq\lq small" sets in $\r$, $E^c$ and $\s^c$ are complements of $E$ and 
$\s$ and $C$ is independent of
$f$. In particular, it follows that if supp$f \subset E$ and supp$\hat f \subset \s$ then $f 
\equiv 0$.  We will use the following definition of the Fourier 
transform: $\hat f(x) = \int f(y) e^{-i2\pi x\cdot y}dy$ and the corresponding inverse Fourier 
transform: $\check f(x) = \int f(y) e^{i2\pi x\cdot y}dy$. 
There are several examples of the Uncertainty 
Principle of 
form \bref{UPIN}. One of them is
the Amrein-Berthier theorem (\cite{HJ}, p.97), \cite{AB} (which is a quantitative version of a 
result due to Benedicks \cite{B}).
In this
theorem sets of finite measure play the role of small sets, i.e., if a function 
$f$ is supported
on a set of finite measure then $\hat f$ can not be concentrated on
a set of finite measure unless $f$ is the zero function. The quantitative 
version of this theorem says that
$$\int |f|^2 \le C(\int\limits_{E^c} |f|^2 + \int\limits_{\Sigma^c}|\hat f|^2)$$
where $E$ and $\Sigma$ are sets of finite measure and the constant $C$ doesn't
depend on $f$. It is interesting to note that the optimal estimate of $C$, which
depends only on measures $|E|$ and $|F|$, was obtained by F. Nazarov relatively
recently \cite{N}.\\

Two more examples of the Uncertainty Principle which are of particular interest to us are
the Logvinenko-Sereda theorem (\cite{HJ}, p.112), \cite{LS} and Wolff's theorem \cite{W}. 
In the case
of the Logvinenko-Sereda theorem  compact sets and
the complements of relatively dense subsets  play the role of 
small sets. A measurable 
set $E^c \subset \r$ 
is called
relatively dense  if there exist a disc $D$ and $\gamma > 0$ such that
\begin{eqnarray}
|E^c \cap (D + x)| \ge \gamma \cdot |D|\label{RELDEN}
\end{eqnarray}
for every $x  \in \r$. It is intuitively clear that $E$ is
a small set in a certain sense. The theorem states that if $\hat f$ is
supported in a compact set $\Sigma$ and a set $E^c$ is relatively dense then
\begin{eqnarray}\int |f|^2 \le C(E, \Sigma)\int\limits_{E^c} |f|^2 
\label{E}\end{eqnarray}
where $C(E, \Sigma)$ depends only on $E$ and $\Sigma$ but doesn't depend on $f$. It is a 
well-known fact that relative density \bref{RELDEN}, 
or \lq\lq thickness", of $E^c$ is also necessary for the  inequality \bref{E} to
hold. See for example (\cite{HJ}, p.113). The Logvinenko-Sereda theorem, which
was motivated by the theory of PDE's, is a result of the theory of entire 
functions of exponential type.
In his earlier paper \cite{K1} the author found a new proof of the Logvinenko-Sereda theorem and 
obtained a sharp 
estimate
of $C$ in \bref{E} which is  polynomial in terms of the density $\gamma$: 
$C \sim \frac 1 {\gamma ^ a}$ rather than a previously known exponential one:  
$C \sim \exp {\frac 1 {\gamma ^ b}}$ (see for example the text book 
by Havin and Joricke (\cite{HJ}, p.112) ). The author showed that his estimate is 
optimal in terms of $\gamma$,  size of $D$ and size of $\Sigma$. 
The author further generalized the inequality \bref{E}
for the case of Fourier transforms supported in the union of finitely many 
compact sets ($\hat f \subset \bigcup\limits_{i = 1}^{n} (\Sigma + \lambda_i)$)
with an estimate of $C$ depending
only on the number of the sets but not how they are placed \cite{K1},
\cite{K3}. In his other paper \cite{K4} the author extended the inequality \bref{E} to 
non-compactly supported Fourier transforms
which are supported in an infinite 
sequence of lacunarily-placed compact sets 
($\hat f \subset \bigcup\limits_{i = -\infty}^{\infty} (\Sigma + \lambda_i)$ where 
$\Lambda = \{\lambda_i \}_{i = - \infty}^{\infty}$  is a lacunary sequence
in $\r$).\\

In Wolff's theorem (\cite{W}, Theorem 2.1)
so called $\epsilon$-thin sets rather than relatively dense ones play the role 
of small sets. Let $\rho(x) = \min {(1, \frac 1 {|x|})}$. A set 
$E \subset \r$ is called $\epsilon$-thin if 
$$|E \cap D(x, \rho(x))| \le \epsilon |D(x, \rho(x))| $$
for all $x \in \r$, where $D(x, r)$ is the disc centered at $x$ with radius $r$. 
The theorem says that if $\epsilon$ is small enough and 
$E$ and $\Sigma$ are $\epsilon$-thin then
\begin{eqnarray}\int |f|^2 \le C(\int\limits_{E^c} |f|^2 + 
\int\limits_{\Sigma^c}|\hat f|^2)\label{UP}\end{eqnarray}
where  $C$ is a universal constant.\\

Our main result is a new version of the Uncertainty Principle which links the Logvinenko-Sereda 
theorem and Wolff's theorem.
Suppose $\rho_1: {\Bbb {R^+}} \rightarrow {\Bbb {R^+}}$ and $\rho_2: {\Bbb {R^+}} \rightarrow 
{\Bbb {R^+}}$ are continuous non-increasing functions and
there exist  $C_1 > 0$, $C_2 > 0$ such that
\begin{eqnarray}
\frac {C_2} {\rho_2 (\frac {C_1} {\rho_1(t)})} \ge t \label{RHO21}
\end{eqnarray}
for all $t \ge 0$. Since the functions $\rho_1$ and $\rho_2$ are continuous 
and non-increasing, we 
also have 
\begin{eqnarray*}
\frac {C_1} {\rho_1 (\frac {C_2} {\rho_2(t)})} \ge t
\end{eqnarray*}
for all $t \ge 0$.
Clearly, $\rho_1 (t) \rightarrow 0$ and  $\rho_2 (t) \rightarrow 
0$ when 
$t \rightarrow \infty$. As an example we can take $\rho_1 (t) = 
\min (\frac 1 {t^a}, 1)$ and $\rho_2(t) = \min(\frac 1 {t^{\frac 1 a}}, 1)$ where $a >0$. Denote 
by 
$D(x, r)$ the disc centered at $x \in \r$ of radius $r$. We call a pair of sets $E, \s \in \r$  
$\epsilon$-thin with respect to the pair of functions $\rho_1$ 
and 
$\rho_2$ if 
\begin{eqnarray}
&&|E \cap D(x, \rho_1(|x|))| \le \epsilon |D(x, \rho_1(|x|))| \nonumber \\
\mbox{and} \label{thin} \\
&&|\s \cap D(x, \rho_2(|x|))| \le \epsilon |D(x, \rho_2(|x|))| \nonumber
\end{eqnarray}
for all $x \in \r$. Note that the complements of $E$ and $\s$ possess some sort of density 
$\gamma = 1 - \epsilon$ with respect to discs $D(x, \rho_1(|x|))$ 
and $D(x, \rho_2(|x|))$ correspondingly. Our main result is the following
\begin{thm}If $\rho_1$ and $\rho_2$ satisfy \bref{RHO21} then there exist $\epsilon > 0$ and $C > 
0$ such that for any pair of $\epsilon$-thin sets $E$ and 
$\s$ with respect to $\rho_1$ and  $\rho_2$ as in \bref{thin} we have
\begin{eqnarray}\int |f|^2 \le C(\int\limits_{E^c} |f|^2 + 
\int\limits_{\Sigma^c}|\hat f|^2) \label{theorem}\end{eqnarray}
for every $f \in L^2$.
\end{thm}
Note that the inequality \bref{RHO21} is scale invariant in the following sense. If we replace 
$f(x)$ with $\tilde f(x) = f(kx)$ then replace $E$ with $\tilde E = kE$ and $\s$ with
$\tilde{\s} = \frac 1 k \s$ where $k > 0$. Then $\tilde E$ and $\tilde{\s}$ are $\epsilon$-thin 
with respect to $\tilde {\rho}_1 (t) = k\rho_1(\frac t k)$ and 
$\tilde {\rho}_2 (t) = \frac 1 k \rho_2(kt)$. It is easy to check that the inequality 
\bref{RHO21} is preserved: $\frac {C_2} {\tilde {\rho}_2 (\frac {C_1} 
{\tilde{\rho}_1(t)})} \ge t$ for all $t \ge 0$.\\

As an application of {\bf Theorem 1} we obtain the following result which is a generalization of 
Theorem 2.3 in \cite{W}. If $G$ is a function on $\r$, 
$\|G\|_{\infty} = 1$, define an operator $T_{G}: L^2 \rightarrow L^2$ as
$$T_G f = \widehat{G f}.$$
Then $\|T_G\|_2 = 1$. The next theorem shows that under certain condition on $G$ and $H$ we have 
$\|T_H T_G\|_2 < 1$.
\begin{thm}
Suppose $\mu_1$ and $\mu_2$ are probability measures on the real line, which are not 
$\delta$-measures. Let $G$ and $H$ be functions on $\r$ satisfying
$$|G(x)| \le |\hat \mu_1(|x|^p)|,$$
$$|H(x)| \le |\hat \mu_2(|x|^{p'})|$$
where $\frac 1 p + \frac 1 {p'} = 1$, $1 < p < \infty$. Then
$$\|T_H T_G\|_2 \le \beta < 1$$
where $\beta$ depends only on $\mu_1$, $\mu_2$ and $p$.
\end{thm}

\section{Proof of Theorems 1 and 2}\label{proof1}
We will construct a pair of bounded operators $S$ and $T$ on $L^2(\r)$   satisfying the following 
conditions:
$\| S \chi_E f\|_2^2 \le \alpha(\epsilon)\|f\|_2^2$ and 
$\|\chi_{\s} \widehat{ T f}\|_2^2 \le 
\beta(\epsilon)\|f\|_2^2$ 
where $\alpha(\epsilon) \lesssim (C_1^d + 1)\epsilon \rightarrow 0$ and $\beta(\epsilon) \lesssim 
(C_2^d + 1)\epsilon \rightarrow 0$ as 
$\epsilon \rightarrow 0$ and such that $S + T$ is the identity operator.  Then
\begin{eqnarray}
\|f\|_2^2 = \|\hat f\|_2^2 &=& \|\chi_{\s^c}\hat f\|_2^2 + \|\chi_{\s}\hat f\|_2^2\nonumber \\
&=&\|\chi_{\s^c}\hat f\|_2^2 + \|\chi_{\s}(\widehat {Sf} + \widehat {Tf})\|_2^2\nonumber \\
&\le& \|\chi_{\s^c}\hat f\|_2^2 + 2\|\chi_{\s}\widehat {Sf}\|_2^2 + 2\|\chi_{\s}\widehat 
{Tf}\|_2^2\nonumber \\
&\le& \|\chi_{\s^c}\hat f\|_2^2 + 2\|Sf\|_2^2 + 2\beta(\epsilon)\|f\|_2^2\nonumber \\
&\le& \|\chi_{\s^c}\hat f\|_2^2 + 4\|S\chi_{E^c}f\|_2^2 + 4\|S\chi_{E}f\|_2^2 + 
2\beta(\epsilon)\|f\|_2^2\nonumber \\
&\le& \|\chi_{\s^c}\hat f\|_2^2 + C \|\chi_{E^c}f\|_2^2 + 4(\alpha(\epsilon) + 
\beta(\epsilon))\|f\|_2^2. \label{lasteps}
\end{eqnarray}
We used here the fact that $S$ is a bounded operator on $L^2$ and therefore $\|S\chi_{E^c}f\|_2 
\lesssim \|\chi_{E^c}f\|_2$.
If $\epsilon$ is small enough so that $4(\alpha(\epsilon) + \beta(\epsilon)) \le \frac 1 2$ then 
we obtain the desired result
\begin{eqnarray*}
\int |f|^2 \lesssim \int\limits_{E^c} |f|^2 +  \int\limits_{\s^c}|\hat f|^2.
\end{eqnarray*}
The theorem contains Wolff's theorem and the Logvinenko-Sereda theorem mentioned in the 
introduction as two extreme cases. To show this
consider $\rho_1 (t) = 
\min (\frac 1 {t^{\frac 1 n}}, 1)$ and $\rho_2(t) = \min(\frac 1 {t^{n}}, 1)$ where $n >0$ then 
$\frac 1 {\rho_2(\frac 1 {\rho_1(t)})} \ge t$. Here $C_1 = C_2 = 1$. If $n = 1$ we 
get Wolff's result. The Logvinenko-Sereda theorem is obtained if we let $n \rightarrow \infty$ 
since $\rho_1(t) \rightarrow 1$ and $\rho_2(t) \rightarrow 
\chi_{[0,1]}(t)$ as $n \rightarrow \infty$, i.e., if
$E$ is the complement of a relatively dense set such that $|E \cap (D + x)| \le \epsilon |D|$ for 
every $x \in \r$ where $D = D(0,1)$ then $E$ is 
$\epsilon$-thin with respect to $\rho_1(t) = 1$ and if a set $\s \subset \epsilon\cdot D$ then 
$\s$ is $\epsilon$-thin with respect to $\rho_2(t) = 
\chi_{[0,1]}(t)$. Note that if sets $E$ and $\s$ are $\epsilon$-thin with respect to $\rho_1(t) = 
1$ and $\rho_2(t) = \chi_{[0,1]}(t)$ correspondingly then 
$E$ is the complement of a relatively dense set and $\s \subset D$. Now we will justify passing 
to the limit. Suppose  $E$ is the complement of a relatively 
dense set such that $|E \cap (D + x)| \le \epsilon |D|$ for every $x \in \r$ and a set $\s 
\subset \epsilon\cdot D$. Let $E_n = E \cap D(0,n)$ and $\s_n = \s$ 
then $E_n$ and $\s_n$ are 
$2\epsilon$-thin with respect to $\rho_1 (t) = 
\min (\frac 1 {t^{\frac 1 n}}, 1)$ and $\rho_2(t) = \min(\frac 1 {t^{n}}, 1)$ correspondingly for 
large enough $n$. Then there exists $\epsilon > 0$ such that 
for all large enough $n$ we have
$$\int |f|^2 \le C( \int\limits_{E_n^c} |f|^2 +  \int\limits_{\s_n^c}|\hat f|^2)$$
where $C$ does not depend on $n$ since all constants in the theorem are uniform in $n$ for these 
pairs of $\rho_1$ and $\rho_2$.
Using $\int\limits_{E_n^c} |f|^2 
\rightarrow \int\limits_{E^c} |f|^2$ and  $\int\limits_{\s _n^c} |f|^2 
= \int\limits_{\s ^c} |\hat f|^2$  we get the Logvinenko-Sereda theorem
\begin{eqnarray*}
\int |f|^2 \lesssim \int\limits_{E^c} |f|^2 +  \int\limits_{\s^c}|\hat f|^2.
\end{eqnarray*}
Now we will show how to construct a pair of such operators $S$ and $T$. We will use some 
technique from (\cite{W}, Theorem 2.1). 
Let $\psi_0: \r \rightarrow {\Bbb {R}}$ be a radial Schwartz function supported in $D(0,2)$ such 
that $\psi_0 \equiv 1$ in $D(0,1)$ and $0 \le \psi_0 \le 1$. 
In addition, we assume that $q(|x|) = \psi_0(x)$ is a non-increasing function of $|x|$ where 
$q(r)$ is a function defined on ${\Bbb {R^+}}$. We will use the 
last property only to prove \bref{Lx} in {\bf Lemma 3}. Define $\psi_j(x) = 
\psi_0(\frac x {2^j}) - \psi_0(\frac x {2^{j-1}})$ for integer $j \ge 1$. It is clear that 
$\psi_j(x)$ is supported in $2^{j-1} \le |x| \le 2^{j+1}$ and 
$\sum\limits_{j \ge 0} \psi_j \equiv 1$. Define $\phi = \check{\psi_0}$ and $\phi_j(x) = C_1^d 
\rho_1^{-d}(2^j)\phi(\frac {C_1x} {\rho_1(2^j)})$. Define the 
operators $S$ and $T$ on $L^2(\r)$ in the following way:
\begin{eqnarray}
S f = \sum\limits_{j \ge 0} \psi_j \cdot (\phi_{j-1} * f) \label{S}
\end{eqnarray}
and
\begin{eqnarray}
T f = \sum\limits_{j \ge 0} \psi_j \cdot (f - \phi_{j-1} * f). \label{T}
\end{eqnarray}
Note that the infinite sums in \bref{S} and \bref{T} converge pointwise since they have at most 
three nonvanishing terms at a given point. It is also clear 
that
$S f + T f \equiv f$.  We have
\begin{eqnarray*}
S f (x) = \int K(x,y) f(y) dy
\end{eqnarray*} 
where 
\begin{eqnarray}
K(x,y) = \sum\limits_{j \ge 0} \psi_j (x)\phi_{j-1}(x-y). \label{K}
\end{eqnarray}
We also have
\begin{eqnarray*}
\widehat {T f} (x) = \int L(x,y) \hat f(y) dy
\end{eqnarray*} 
where 
\begin{eqnarray}
L(x,y) &=& \sum\limits_{j \ge 0} \hat \psi_j (x - y) (1 - \hat\phi_{j-1}(y)) \label{L1} \\
&=& \phi (x - y) (1 - \hat\phi_{-1}(y)) + \sum\limits_{j \ge 1} (2^{jd}\phi(2^j(x - y)) - 
2^{(j-1)d}\phi(2^{j-1}(x - y)))(1 - \hat\phi_{j-1}(y))\nonumber \\
&=& \sum\limits_{j \ge 0} 2^{jd}\phi(2^j(x - y))(\hat\phi_{j}(y) - \hat\phi_{j-1}(y))\nonumber \\
&=& \sum\limits_{j \ge 0} 2^{jd}\phi(2^j(x - y))(\psi_0({\rho_1(2^{j})}y/C_1) - 
\psi_0({\rho_1(2^{j-1})}y/C_1)). \label{L}
\end{eqnarray}
We used here $\hat \psi_j(z) = 2^{jd}\hat\psi_0(2^j z) - 2^{(j-1)d}\hat\psi_0(2^{j - 1} z) = 
2^{jd}\phi(2^j z) - 2^{(j-1)d}\phi(2^{j-1}z)$ for $j \ge 1$, 
summation by parts and $\hat\phi_{j}(y) = \psi_0({\rho_1(2^{j})}y/C_1)$.
Note that for a fixed $y$ the sums in \bref{L1} and \bref{L} have only finitely many terms since 
$\hat \phi_j(y) = \psi_0({\rho_1(2^{j})}y/C_1) = 1$ if $|y| 
\le \frac {C_1} {\rho_1(2^j)}$ and $\frac {C_1} {\rho_1(2^j)} \rightarrow \infty$ when $j 
\rightarrow \infty$. The four lemmas below are analogous to Lemma 
2.2 in \cite{W}.\\
Now we will show that $S$ is a bounded operator on $L^2$.
It will suffice to prove the following lemma: 
\begin{lem}
\begin{eqnarray}
\sup\limits_x\int|K(x,y)|dy \le C \label{Kx}
\end{eqnarray} 
and 
\begin{eqnarray}
\sup\limits_y\int|K(x,y)|dx \le C \label{Ky}.
\end{eqnarray}
where $C$ is an absolute constant which does not depend on $\rho_1$.
\end{lem}
{\bf Proof of Lemma 1.}
\bref{Kx} follows from the facts that for a fixed $x$ the sum in \bref{K} contains at most three 
nonvanishing terms, $|\psi_j| \le 1$ and $\|\phi_j\|_1 = 
\|\phi\|$. 
Therefore, 
\begin{eqnarray*}
\sup\limits_x\int|K(x,y)|dy \le 3\|\phi\|_1. 
\end{eqnarray*}
Fix $y$ and note that there are at most four values of $j$ such that dist$(y,\mbox{supp }\psi_j) 
< 2^{j-2}$. Call this set of 
$j$'s A. We have
\begin{eqnarray}
\int|K(x,y)|dx &\le& 4\|\phi\|_1 + \sum\limits_{j \notin A}\int|\psi_j (x)| \cdot 
|\phi_{j-1}(x-y)| dx\nonumber \\
&\le& 4\|\phi\|_1 + \sum\limits_{j \notin A}\int|\psi_j (x)| \cdot \frac 
{C\rho_1^{-d}(2^{j-1})C_1^d} {(1 + \frac {C_1|x-y|}{\rho_1(2^{j-1})})^{2d}} 
dx\nonumber \\
&\le& 4\|\phi\|_1 + \sum\limits_{j \notin A}\int|\psi_j (x)| \cdot \frac 
{C\rho_1^{-d}(2^{j-1})C_1^d} {(1 + \frac {C_12^{j-2}}{\rho_1(2^{j-1})})^{2d}} 
dx\nonumber \\
&\le& 4\|\phi\|_1 + C\sum\limits_{j \ge 0}2^{jd} \frac {\rho_1^{-d}(2^{j-1})C_1^d} {(1 + \frac 
{C_12^{j-2}}{\rho_1(2^{j-1})})^{2d}}\label{Ky2}. 
\end{eqnarray}
To estimate the second term in \bref{Ky2} we will use the fact that $\rho_1$ is non-increasing. 
 Choose the smallest integer $k \ge -1$ such that $C_12^{k} 
\ge \rho_1(2^{k})$ and split the sum  correspondingly into two parts (ignore the first part if $k 
= -1$)
\begin{eqnarray*}&&\sum\limits_{j \ge 0}^{k}2^{jd} \frac {\rho_1^{-d}(2^{j-1})C_1^d} {(1 + \frac 
{C_12^{j-2}}{\rho_1(2^{j-1})})^{2d}} + \sum\limits_{j \ge 
k + 1}2^{jd} \frac 
{\rho_1^{-d}(2^{j-1})C_1^d} {(1 + \frac {C_12^{j-2}}{\rho_1(2^{j-1})})^{2d}} \le \\
&&C 2^{kd}\rho_1^{-d}(2^{k-1})C_1^d + C\sum\limits_{j \ge k + 1}2^{-jd}C_1^{-d}\rho_1^{d}(2^{k}) 
\le C.
\end{eqnarray*}
Therefore, it follows from \bref{Ky2} that 
$$\sup\limits_y\int|K(x,y)|dx \le C$$ 
where $C$ does not depend on $\rho_1$. \hfill$\square$\\

Thus, $S$ is a bounded operator on $L^2$ whose norm does not depend on $\rho_1$. \\
Now we will show that 
\begin{eqnarray}
\|S \chi_E f\|_2 \le C \sqrt {\epsilon}\|f\|_2. \label{SE}\end{eqnarray}
Since we have already shown that
\begin{eqnarray*}
\sup\limits_y\int|K(x,y)|dx \le C. 
\end{eqnarray*}
it will suffice to prove the next lemma.
\begin{lem} 
\begin{eqnarray}
\sup\limits_x\int\limits_E |K(x,y)|dy \le C \epsilon\label{Keps}
\end{eqnarray}
where $C \lesssim C_1^d + 1$.
\end{lem}
{\bf Proof of Lemma 2.}
To obtain \bref{Keps} we will need the following geometrical property:
\begin{eqnarray}
|D(x, r) \cap E| \le C\epsilon |D(x,r)| \label{DE}
\end{eqnarray}
for all $x$ and $r \ge \rho_1(|x|)$.
This inequality is based on the fact that we can cover $\bar D(x,r)$ by disks $D(x_i, \rho_1(|x_i|))$, 
$x_i \in \bar D(x,r)$, such that
\begin{eqnarray}
\sum\limits_i |D(x_i, \rho_1(|x_i|))| \le C |D(x,r)|. \label{Dxi}
\end{eqnarray}
To show this we will use only the continuity of $\rho_1$. First we claim that $\bar D(x, r)$ can be covered 
by balls $D(x_k, \rho_1(|x_k|)/3)$ with $x_k \in \bar D(x,r)$ and $\rho_1(|x_k|) \le 
3r$. 
Suppose towards a contradiction that there is $y \in \bar D(x,r)$ which is not covered. Then 
$\rho_1(|y|) > 3r$. Consider $h(t) = \rho_1(|(1-t)x + ty|)$, $0 \le t 
\le 1$. It is continuous and $h(0) \le r$ and $h(1) > 3r$. Pick $t$ such that $h(t) = 3r$. Let $z 
= (1-t)x + ty$ then $z \in D(x,r)$,  
$\rho_1(|z|) = 3r$ and $y \in D(z, \rho_1(|z|) / 3)$. It gives us a contradiction. Thus, $\bar D(x, r)$ can be covered 
by balls $D(x_k, \rho_1(|x_k|)/3)$ with $x_k \in \bar D(x,r)$ and $\rho_1(|x_k|) \le 
3r$. 
Next we choose a finite subcover of $\bar D(x, r)$ (a compact set) by these balls. Applying a well-known result on 
covering (see 
(\cite{R}, 7.3)) we can choose a disjoint subcollection of these balls $D(x_i, \rho_1(|x_i|)/3)$ such 
that $D(x,r) \subset \bigcup\limits_i D(x_i, \rho_1(|x_i|))$. 
We 
also have $\bigcup\limits_i D(x_i, 
\rho_1(|x_i|)/3) \subset D(x, 2r)$. Therefore, 
$\sum\limits_i |D(x_i, \rho_1(|x_i|))| = 3^d\sum\limits_i |D(x_i, \rho_1(|x_i|)/3)| = 3^d 
|\bigcup\limits_i D(x_i, \rho_1(|x_i|)/3)| \le 
 3^d|D(x,2r)| = C|D(x,r)|$ which gives us \bref{Dxi}. Now \bref{DE} follows from \bref{Dxi}:
 $$|D(x, r) \cap E| \le \sum\limits_i |D(x_i, \rho_1(|x_i|)) \cap E| \le \epsilon  \sum\limits_i 
|D(x_i, \rho_1(|x_i|))| \le C\epsilon |D(x,r|.$$
 Since $K(x,y) = \sum\limits_{j \ge 0} \psi_j (x)\phi_{j-1}(x-y)$ has at most three nonvanishing 
terms for each fixed $x$ and $|\psi_j(x)| \le 1$, to prove 
\bref{Keps} it is enough to 
show that
\begin{eqnarray}
\sup\limits_x\int\limits_E |\phi_j(x-y)|dy \le C \epsilon\label{Keps1}
\end{eqnarray}
when $x \in$ supp $\psi_{j+1}$, i.e., $2^{j} \le |x| \le 2^{j+2}$ and therefore $\rho_1(|x|) \le 
\rho_1(2^j)$. Let $k_1 \ge 0$ be the 
smallest integer such that $2^{k_1} \ge C_1$. Therefore, using 
\bref{DE} we get \bref{Keps1}:
\begin{eqnarray*}
\int\limits_E |\phi_j(x-y)|dy &=& \int\limits_E |\rho_1^{-d}(2^j)C_1^d\phi(\frac {C_1(x - y)} 
{\rho_1(2^j)})|dy\\
&\le& C\int\limits_E \frac {\rho_1^{-d}(2^j)C_1^d}{(1 + \frac {C_1|x - y|} 
{\rho_1(2^j)})^{2d}}dy\\ 
&\lesssim& \int\limits_{E \cap D(x, 2^{k_1} \rho_1(|x|)/C_1)} + \sum\limits_{k > 
k_1}\int\limits_{E \cap D(x, 2^k \rho_1(|x|)/C_1)\backslash D(x, 2^{k-1} 
\rho_1(|x|)/C_1)}\frac {\rho_1^{-d}(2^j)C_1^d} {(1 + \frac {C_1|x - y|} {\rho_1(2^j)})^{2d}}dy\\
&\le& \epsilon (2^{k_1} \rho_1(|x|)/C_1)^d{\rho_1^{-d}(2^j)C_1^d} + \sum\limits_{k > k_1}\epsilon 
(2^k \rho_1(|x|)/C_1)^d \frac {\rho_1^{-d}(2^j)C_1^d} {(1 + 
\frac {2^k\rho_1(|x|)}{\rho_1(2^j)})^{2d}}\\
&\lesssim& \epsilon (C_1^d + 1). 
\end{eqnarray*} \hfill$\square$\\
Thus, $\|S \chi_E f\|_2^2 \lesssim (C_1^d + 1) {\epsilon}\|f\|_2^2.$\\
Since $S + T = I$, it follows that  $T$ is also a bounded operator on $L^2$. 
However, we will need the following lemma for the operator $T$ later which is analogous for $\bf 
Lemma 1$ for the operator $S$:
\begin{lem} 
\begin{eqnarray}
\sup\limits_y\int|L(x,y)|dx \le C \label{Ly}.
\end{eqnarray}
and
\begin{eqnarray}
\sup\limits_x\int|L(x,y)|dy \le C \label{Lx}.
\end{eqnarray}
where $C$ is an absolute constant and does not depend on $\rho_1$.
\end{lem}
{\bf Proof of Lemma 3.}
Recall that $\psi_0(x) = q(|x|)$ where $q(r)$ is a function defined on ${\Bbb {R^+}}$. 
Using \bref{L} we have
\begin{eqnarray*}
\int|L(x,y)|dx &\le& \|\phi\|_1 \sum\limits_{j \ge 0}|\psi_0({\rho_1(2^{j})}y/C_1) - 
\psi_0({\rho_1(2^{j-1})}y/C_1)| \\ 
&\le& C \sum\limits_{j \ge 0} \int\limits_{ {\rho_1(2^{j})}|y|/C_1}^{{\rho_1(2^{j-1})}|y|/C_1} 
|q'(t)| dt \\
&\le& C \int\limits_{0}^{\infty} |q'(t)| dt \le C.
\end{eqnarray*}
We used here the facts that $\psi_0$ is a radial Schwartz function and $\rho_1$ is 
non-increasing.\\
Now we will show that 
\begin{eqnarray*}
\sup\limits_x\int|L(x,y)|dy \le C.
\end{eqnarray*}
We have that $\phi = 
\check{\psi_0}$ is a real-valued radial Schwartz function. Actually, $\phi$ can be extended to an 
entire function of exponential type on ${\Bbb{C}}^d$. It is 
clear that $\widehat{|\phi|}(x)$ is a radial function. Here we took the Fourier transform of the 
absolute value of $\phi$. Denote $p(|x|) = 
\widehat{|\phi|}(x)$ where $p(r)$ is a function defined on ${\Bbb {R^+}}$. Recall from our 
definition of $\psi_0(x)$ that $q(|x|) = \psi_0(x)$ is a 
non-increasing function of 
$|x|$. Therefore, $0 \le \psi_0({\rho_1(2^{j})}y/C_1) - \psi_0({\rho_1(2^{j-1})}y/C_1) = 
\hat\phi_{j}(y) - \hat\phi_{j-1}(y)$ since $\rho_1$ is a 
non-increasing function too.
Let 
$$
L_n(x,y) = \sum\limits_{j = 0}^{n} 2^{jd}\phi(2^j(x - y))(\psi_0({\rho_1(2^{j})}y/C_1) - 
\psi_0({\rho_1(2^{j-1})}y/C_1)).
$$
Since $L_n(x,y) \rightarrow L(x, y)$, it is enough to show that $\int |L_n(x,y)| dy \le C$ to get 
\bref{Lx} by an application of Fatou's Lemma. We have
\begin{eqnarray}
\int |L_n(x,y)| dy &\le& \int \sum\limits_{j = 0}^{n} 2^{jd}|\phi(2^j(x - y))| \cdot 
(\hat\phi_{j}(y) - \hat\phi_{j-1}(y)) dy\nonumber \\
&=& \int \sum\limits_{j = 0}^{n} 2^{jd}{|\phi(2^j(x - \cdot))|}^{\wedge}(z) \cdot (\phi_{j}(z) - 
\phi_{j-1}(z)) dz\nonumber \\
&=& \int \sum\limits_{j = 0}^{n} e^{-i2\pi x z}\widehat{|\phi|}(\frac z {2^j}) \cdot (\phi_{j}(z) 
- \phi_{j-1}(z)) dz\nonumber \\
&=& \int e^{-i2\pi x z}\left(\widehat{|\phi|}(\frac z {2^n}) \phi_{n}(z) - \widehat{|\phi|}(\frac 
z {2^0}) \phi_{-1}(z)\right)dz + \nonumber \\
&& \int\sum\limits_{j = 1}^{n} 
e^{-i2\pi x z}(\widehat{|\phi|}(\frac z {2^{j-1}}) - \widehat{|\phi|}(\frac z {2^j})) \cdot \phi_{j-1}(z) 
dz\nonumber \\
&\le& 2\|\phi\|_1 \cdot \|\psi_0\|_{\infty} + \sum\limits_{j = 1}^{n} \int |\phi_{j-1}(z)| \int\limits_{\frac {|z|} 
{2^j}}^{\frac {|z|} {2^{j-1}}} |p'(t)|dt 
 dz. \label{Lx3}
\end{eqnarray}
Here we used the identity $\int f \hat g = \int \hat f g$ and summation by parts.
Change the order of integration to estimate the second term in \bref{Lx3} by
\begin{eqnarray}
&&\int\limits_0^{\infty} |p'(t)| \sum\limits_{j = 1}^{n} \int\limits_{2^{j-1}t \le |z| \le 2^j t} 
\frac {C_1^d}{\rho_1^d(2^{j-1})}|\phi(\frac 
{C_1 
z}{\rho_1(2^{j-1})})| dz dt\nonumber \\
&=& \int\limits_0^{\infty} |p'(t)| \sum\limits_{j = 1}^{n} \int\limits_{C_1 
2^{j-1}t\rho_1^{-1}(2^{j-1}) \le |z| \le C_1 2^j 
t\rho_1^{-1}(2^{j-1})} 
|\phi(z)| 
dz dt\nonumber \\
&\le& \int\limits_0^{\infty} |p'(t)| \cdot \|\phi\|_1 dt.\label{Lx4}
\end{eqnarray}
Here we used the fact that $\rho_1$ is non-increasing.
All we need to show now is that $\int\limits_0^{\infty} |p'(t)|dt < \infty$. To obtain this we 
will prove that $|p'(t)| 
\lesssim t^{-2}$ (actually $\lesssim t^{-2 - \frac{d-1}{2}}$.) Recall that $\phi$ is a 
real-valued radial Schwartz function which can be extended to an entire 
function of exponential type on ${\Bbb{C}^d}$. Let $g(t) = \phi(t,0,0,...,0)$ then $g(t)$ is an 
even real-valued Schwartz function on ${\Bbb{R}}$ which can be 
extended to an entire 
function of exponential type on ${\Bbb{C}}$ and $g(|x|) = \phi(x)$.
Let 
$\alpha_k$, $k = 1, 2, ...$, be the positive roots of $g$ in increasing order where the function 
$g$ changes its sign. Let $\alpha_0 = 0$ (note $g(0) = 
\phi(0,0,...,0) > 1$). We have
\begin{eqnarray}
p(t) = \sum\limits_{k \ge 0}(-1)^k\int\limits_{\alpha_k}^{\alpha_{k+1}} g(s) 
\hat{d\sigma}(st)s^{d-1}ds \label{hatphi}
\end{eqnarray}
where 
\begin{eqnarray}
\hat{d\sigma}(|x|) = \int\limits_{|\xi| = 1}e^{-i2\pi x\cdot \xi}d\sigma(\xi) = 
|\sigma_{d-2}|\cdot \int\limits_{0}^{\pi}e^{-i2\pi |x| \cos {\theta}}(\sin 
{\theta})^{d-2}d \theta
\label{hatdsig}
\end{eqnarray}
is a real-valued radial function. Here $|\sigma_{d-2}|$ is the area of the $(d-2)$-dimensional 
sphere. It is well-known that 
$\hat{d\sigma}(r) = Re(B(r))$ where $B(r) = a(r)e^{i2\pi r}$ and $a(r)$ satisfies the following 
estimates:
$$|a^{(k)}(r)| \le \frac C {r^{\frac {d-1} {2} + k}}.$$
See, for example, (\cite{S}, Ch. viii).\\
Denote by $$F(r) = \int\limits_0^{r}\hat{d\sigma}(t) dt.$$
We need only the facts that $|\hat{d\sigma}(r)| \le C$ and $|F(r)| = 
|\int\limits_0^{r}\hat{d\sigma}(t) dt| \le C$ which easily follow from the definition 
\bref{hatdsig} and integration by parts (if $d = 1$ then $\hat{d\sigma}(r) = 
\cos{2\pi r}$ 
and the estimate is obvious). Actually, it follows from the properties of $\hat{d\sigma}$ that
$|\int\limits_r^{\infty}\hat{d\sigma}(t) dt| \lesssim r^{-\frac{d-1}{2}}$.) Differentiating 
\bref{hatphi} we get
\begin{eqnarray}
p'(t) &=& \sum\limits_{k \ge 0}(-1)^k\int\limits_{\alpha_k}^{\alpha_{k+1}} g(s) 
(\hat{d\sigma})'(st)s^{d}ds = \nonumber \\
&& \sum\limits_{k \ge 0}-(-1)^k\frac 1 t \int\limits_{\alpha_k}^{\alpha_{k+1}} (s^{d}g(s))' 
\hat{d\sigma}(st)ds = \nonumber \\
&& \sum\limits_{k \ge 0}(-1)^{k+1}\frac 1 {t^2} \int\limits_{\alpha_k}^{\alpha_{k+1}} 
(s^{d}g(s))' d(\int_0^{st}\hat{d\sigma}(r)dr) = \nonumber \\
&& \sum\limits_{k \ge 0}(-1)^{k+1}\frac 1 {t^2}\left ( \left.(s^d g(s))'F(st)\right \vert_{s = 
\alpha_k}^{s = \alpha_{k+1}}  - 
\int\limits_{\alpha_k}^{\alpha_{k+1}} (s^d g(s))''F(st)ds\right ).\nonumber
\end{eqnarray}
We integrated by parts twice and used that $\alpha_0 = 0$ and $g(\alpha_k) = 0$ for $k = 1, 2, 
...$.
Since $g$ is an entire function of exponential type, we can bound the number of its roots in the 
interval $[0, 2^j]$ by $\lesssim 2^j$ which follows from 
Jensen's formula (\cite{R}, p. 309). Thus,
\begin{eqnarray}
|p'(t)| &\lesssim& \frac 1 {t^2}\sum\limits_{\alpha_k < 1} |(s^d g(s))'|(\alpha_k) + \frac 1 
{t^2}\sum\limits_{j \ge 0} \sum\limits_{2^j \le 
\alpha_k < 2^{j+1}} |(s^d g(s))'|(\alpha_k) + \frac 1 {t^2}\int |(s^d g(s))''|ds\nonumber \\
&\lesssim& \frac 1 {t^2} + \frac 1 {t^2}\sum\limits_{j \ge 0} \frac {2^{jd}2^j} {1 + 2^{2j(d+1)}} 
+ \frac 1 {t^2}\nonumber \\
&\lesssim& \frac 1 {t^2}. \label{Dhatphi3}
\end{eqnarray}
Therefore,
\begin{eqnarray}
\int |p'(t)|dt < \infty.\label{Dhatphi2}
\end{eqnarray}
Substituting \bref{Dhatphi2} into \bref{Lx4} and using \bref{Lx3} we obtain
$$\int |L_n(x,y)| dy \le C$$
which gives us the desired estimate \bref{Lx}:
$$\int |L(x,y)| dy \le C.$$
If we use $|\int\limits_r^{\infty}\hat{d\sigma}(t) dt| \lesssim r^{-\frac{d-1}{2}}$ then we can 
improve \bref{Dhatphi3} 
$$|p'(t)| \lesssim \frac 1 {t^{2 + \frac{d-1}{2}}}.$$
\hfill$\square$\\

In particular, it follows that $T$ is a bounded operator on $L^2$ whose norm does not depend on 
$\rho_1$. \\
The proof of 
\begin{eqnarray}
\|\chi_{\s}\widehat{T f}\|_2 \le C \sqrt {\epsilon}\|f\|_2 \label{Tsigma}\end{eqnarray}
is quite similar to the proof of \bref{SE}.
Since we have already shown that
\begin{eqnarray*}
\sup\limits_x\int|L(x,y)|dy \le C, 
\end{eqnarray*} 
it will suffice to prove the following lemma:
\begin{lem} 
\begin{eqnarray}
\sup\limits_y\int\limits_{\s} |L(x,y)|dx \le C \epsilon\label{Leps}
\end{eqnarray}
where $C \lesssim C_2^d + 1$.
\end{lem}
{\bf Proof of Lemma 4.}
Combining with the argument in the proof of  \bref{Ly} it is enough to show that
\begin{eqnarray}
\int\limits_{\s}2^{jd}|\phi(2^j(x - y))| dx \le C \epsilon \label{Leps1}
\end{eqnarray}
when $\frac {C_1} {\rho_1(2^{j-1})} \le |y| \le \frac {2C_1} {\rho_1(2^{j})}$. Therefore, 
applying \bref{RHO21} we get $\rho_2(|y|) \le \rho_2(\frac {C_1} 
{\rho_1(2^{j-1)}}) \le C_2 2^{-(j-1)}$. Repeating an argument similar to the one to obtain 
\bref{Keps1} we have
\begin{eqnarray*}
\int\limits_{\s}2^{jd}|\phi(2^j(x - y))| dx &\lesssim& \int\limits_{\s \cap D(y, \rho_2(|y|))} + 
\sum\limits_{k > 0}\int\limits_{\s \cap D(y, 2^k\rho_2(|y|)) 
\backslash D(y, 2^{k-1}\rho_2(|y|))} \frac {2^{jd}}{(1 + 2^j|x - y|)^{2d}}dx \\
&\lesssim& 2^{jd}\epsilon (\rho_2(|y|))^d + \sum\limits_{k > 0}\frac {2^{jd}}{(1 + 2^j 2^k 
\rho_2(|y|))^2d}\\
&\lesssim& \epsilon ((2^{j}\rho_2(|y|))^d + 1) \lesssim \epsilon (C_2^d + 1).
\end{eqnarray*}\hfill$\square$\\

Thus, $\|\chi_{\s}\widehat{T f}\|_2^2 \lesssim (C_2^d + 1) {\epsilon}\|f\|_2^2.$\\
Therefore the pair of operators $S$ and $T$ possesses the properties we gave at the beginning of 
the proof of the theorem and we obtain the desired result
\begin{eqnarray*}
\int |f|^2 \lesssim \int\limits_{E^c} |f|^2 +  \int\limits_{\s^c}|\hat f|^2.
\end{eqnarray*}
End of proof of {\bf Theorem 1}. \hfill$\square$\\

It is an open question whether {\bf Theorem 1} holds for every $0 < \epsilon < 1$ but not just 
for sufficiently small $\epsilon$.\\

We will give only a sketch of the proof of {\bf Theorem 2} since it is analogous to the proof of 
Theorem 2.3 in \cite{W}. If $\mu$ is a probability measure 
which 
is not a unit point mass then the set $F = \{|\hat \mu (\xi)| > 1 - \delta\}$ is a complement of 
a relatively dense one (\bref{RELDEN}) if $\delta$ is small 
enough. More precisely, 
$$|F \cap [x, x + 1]| \le \epsilon(\delta)$$ for every $x \in {\Bbb{R}}$, where $\epsilon(\delta) 
\rightarrow 0$ as $\delta \rightarrow 0$. The statement is 
also 
true in $\r$. See for example Lemma 1.1 in \cite{W}. Define similarly $F_1$ and $F_2$ 
correspondingly for $\mu_1$ and $\mu_2$. If $Q_1(x) = |x|^p$ and $Q_2(x) 
= 
|x|^{p'}$ then the sets $E = \{|G(x)| > 1 - \delta\} \subset Q_1^{-1}(F_1)$ and $\s = \{|H(x)| > 
1 - \delta\} \subset Q_2^{-1}(F_2)$ are $\epsilon$-thin with 
respect to $\rho_1(|x|) = \min(\frac 1 {|x|^{p-1}}, 1)$ and $\rho_2(|x|) = \min(\frac 1 {|x|^{p' 
- 1}}, 1)$ correspondingly. Note that $(p-1)(p' - 1) = 1$. 
This 
is true because $Q_1$ and $Q_2$ map discs $D_1 = D(x, \rho_1(|x|)$ and $D_2 = D(x, \rho_2(|x|))$ 
correspondingly onto intervals $I_1(x)$ and $I_2(x)$ of 
length 
$\gtrsim C^{-1}$ since $|\triangledown Q_1 (x)| \sim \frac 1 {\rho_1(|x|)}$ and $|\triangledown 
Q_2 (x)| \sim \frac 1 {\rho_2(|x|)}$ if $|x|$ is large. 
Therefore,
$$\frac {|Q_1^{-1}(F_1) \cap D(x, \rho_1{|x|})|}{|D(x, \rho_1{|x|})|} \lesssim \epsilon^{\frac 1 
p},$$
$$\frac {|Q_2^{-1}(F_2) \cap D(x, \rho_2{|x|})|}{|D(x, \rho_2{|x|})|} \lesssim \epsilon^{\frac 1 
{p'}}$$
since for large $|x|$ 
$$\frac {|Q_1^{-1}(F_1) \cap D(x, \rho_1{|x|})|}{|D(x, \rho_1{|x|})|} \le \frac {|F_1 \cap 
I_1|}{|I_1|}\cdot \frac {\max\limits_{D_1} \frac 1 
{|\triangledown Q_1|}}{\min\limits_{D_1} \frac 1 {|\triangledown Q_1|}} \lesssim \epsilon$$ and 
for small $|x|$ it is $\lesssim \max(\epsilon^{\frac d p}, 
\epsilon)$. Similar estimates 
hold for $Q_2^{-1}(F_2)$.\\

We have 
\begin{eqnarray}
\|Gf\|_2^2 &=& \|Gf\|_{L^2{(E^c)}}^2 + \|Gf\|_{L^2{(E)}}^2\nonumber \\ 
&\le& (1 - \delta)^2\|f\|_{L^2(E^c)}^2 + \|f\|_{L^2(E)}^2\nonumber \\
&=& \|f\|_2^2 - (1 - (1 - \delta)^2)\|f\|_{L^2(E^c)}^2. \label{Gf}
\end{eqnarray}
In a similar way we have 
\begin{eqnarray}
\|H\widehat{Gf}\|_2^2 \le \|\widehat{Gf}\|_2^2 - (1 - (1 - 
\delta)^2)\|\widehat{Gf}\|_{L^2(\s^c)}^2. \label{HGf}
\end{eqnarray}  
It follows from {\bf Theorem 1} that $$\|\widehat{Gf}\|_{L^2(\s^c)}^2 \ge C^{-1}\|Gf\|_2^2 - 
\|Gf\|_{L^2(E^c)}^2$$ where $C \ge 1$. Therefore, \bref{HGf} 
$\lesssim$ 
\begin{eqnarray*}
\|H\widehat{Gf}\|_2^2 &\le& \|Gf\|_2^2 - (1 - (1 - \delta)^2)(C^{-1}\|Gf\|_2^2 - 
\|Gf\|_{L^2(E^c)}^2)\\
&\le& (1 - C^{-1}(1 - (1 - \delta)^2))\|Gf\|_2^2 + (1 - (1 - \delta)^2)(1 - 
\delta)^2\|f\|_{L^2(E^c)}^2\\
&\le& (1 - C^{-1}(1 - (1 - \delta)^2))\|f\|_2^2.
\end{eqnarray*}
We used here \bref{Gf} and the fact  $(1 - (1 - \delta)^2)(1 - C^{-1}(1 - (1 - \delta)^2)) \ge (1 
- (1 - \delta)^2)(1 - \delta)^2$ since $C \ge 1$.
It is clear that $\beta^2 \le 1 - C^{-1}(1 - (1 - \delta)^2) < 1$. 
Instead of $Q_1(x) = |x|^p$ and $Q_2(x) = |x|^{p'}$ we can also use $Q_1(x) = \sum\limits_{i = 
1}^{d} a_i|x_i|^p$ and $Q_2(x) = \sum\limits_{i = 1}^{d} 
b_i|x_i|^{p'}$ where $a_i$ and $b_j$ are non-zero real numbers.
\begin{rem}
If $Q_1(t)$ and $Q_2(t)$ are increasing convex functions satisfying $C_2Q_2'(C_1Q_1'(t)) \ge t$ 
for all $t \ge 0$ and such that $Q_i'(t + \frac 1 {Q_i'(t)}) 
\le C 
Q_i'(t - \frac 1 {Q_i'(t)})$, $i = 1, 2$ for all large $t$ then Theorem 2 holds with $|G(x) \le 
|\hat \mu_1(Q_1(|x|))|$ and $|H(x) \le |\hat 
\mu_2(Q_2(|x|))|$. 
\end{rem}

\section{Counterexamples}\label{counter}
Recall that a pair of sets $E \subset \r$ and $\s \subset \r$ is called $\epsilon$-thin with 
respect to the pair of functions $\rho_1$ and $\rho_2$ 
correspondingly if
\begin{eqnarray*}
|E \cap D(x, \rho_1(|x|))| \le \epsilon |D(x, \rho_1(|x|))|
\end{eqnarray*}
and
\begin{eqnarray*}
|\s \cap D(x, \rho_2(|x|))| \le \epsilon |D(x, \rho_2(|x|))|
\end{eqnarray*}
for all $x \in \r$.
The next lemma shows that the condition \bref{RHO21} $\frac {C_2} {\rho_2(\frac {C_1} 
{\rho_1(t)})} \ge t$ in the {\bf Theorem} is not only sufficient but 
also necessary 
for an inequality of the form
\begin{eqnarray*}
\int |f|^2 \le C(E, \s)(\int\limits_{E^c} |f|^2 + 
\int\limits_{\Sigma^c}|\hat f|^2)
\end{eqnarray*}
to hold for every $f \in L^2$ and every pair of $\epsilon$-thin sets $E$ and $\s$ with respect to 
$\rho_1$ and $\rho_2$ correspondingly.

\begin{lem}
Let $\rho_1: {\Bbb {R^+}} \rightarrow {\Bbb {R^+}}$ and $\rho_2: {\Bbb {R^+}} \rightarrow {\Bbb 
{R^+}}$ be continuous non-increasing functions. Suppose that 
there exists $0 < \epsilon < 1$ such that for every pair of $\epsilon$-thin sets $E$ and $\s$ 
with respect to $\rho_1$ and $\rho_2$ correspondingly we have 
\begin{eqnarray}\int |f|^2 \le C(E, \s)(\int\limits_{E^c} |f|^2 + 
\int\limits_{\Sigma^c}|\hat f|^2)\label{UP1}\end{eqnarray}
for every $f \in L^2$
then there exist $C_2 > 0$ and $C_1 > 0$ such that $\frac {C_2} {\rho_2(\frac {C_1} {\rho_1(t)})} 
\ge t$ for every $t \ge 0$.
\end{lem}

{\bf Proof of Lemma 5:} The pair of functions $\rho_1$ and $\rho_2$ will be fixed throughout the 
proof. First we will show that there exists a universal 
constant $C > 0$ such that \bref{UP1} holds for thin enough sets $E$ and $\s$, i.e., there exist 
$\epsilon_0 \in (0, \epsilon]$ and $C > 0$ such that
\begin{eqnarray*}
\int |f|^2 \le C(\int\limits_{E^c} |f|^2 + 
\int\limits_{\Sigma^c}|\hat f|^2)
\end{eqnarray*}
holds for every pair of thin sets $E$ and $\s$ with thinness $\epsilon_0$ and every $f \in L^2$. 
Suppose towards a contradiction that this is not true. Then 
there exists a sequence of $f_n$ and corresponding $\epsilon_n$-thin sets $E_n$ and $\s_n$ with 
respect to $\rho_1$ and $\rho_2$ with thinness $\epsilon_n < 
\frac {\epsilon} {2^n}$ such that
\begin{eqnarray*}
\int |f_n|^2 > n(\int\limits_{E_n^c} |f_n|^2 + 
\int\limits_{\Sigma_n^c}|\hat f_n|^2).
\end{eqnarray*}
Let $E = \bigcup\limits_{n = 1}^{\infty}E_n$ and $\s = \bigcup\limits_{n = 1}^{\infty}\s_n$. Then 
$E$ and $\s$ are $\epsilon$-thin with respect to $\rho_1$ 
and 
$\rho_2$ correspondingly. On the other hand we have
\begin{eqnarray*}
\int |f_n|^2 > n(\int\limits_{E_n^c} |f_n|^2 + 
\int\limits_{\Sigma_n^c}|\hat f_n|^2) \ge n(\int\limits_{E^c} |f_n|^2 + 
\int\limits_{\Sigma^c}|\hat f_n|^2)
\end{eqnarray*}
which contradicts to \bref{UP1}.\\

Suppose towards a contradiction that for any choice of $C_2 > 0$ and $C_1 > 0$ there exists $t 
\ge 0$ with the property $\frac {C_2} {\rho_2(\frac {C_1} 
{\rho_1(t)})} < t$. 
Then for 
any arbitraly small $\epsilon > 0$ we will construct a sequence of $\epsilon$-thin sets $E_n$ and 
$\s_n$ with respect to $\rho_1$ and $\rho_2$ and a sequence 
of Schwartz functions $f_n$ with  supp $f_n \subset E_n$ such that
\begin{eqnarray*}
\lim\limits_{n \rightarrow \infty} \frac {\int\limits_{\Sigma^c_n}|\hat f_n|^2} {\int |f_n|^2} = 
0. 
\end{eqnarray*}
First we will discuss the 1-dimensional case. To simplify the proof we can assume without loss of 
generality  that $\rho_1 < \frac 1 2$ and $\rho_2 < \frac 1 
2$. Let $\phi$ be a Schwartz function supported in $[-1,1]$. Choose an integer $n > 0$ which we 
will specify later. 
Let $0 < \epsilon < 1$. Note that the functions $\phi(\frac {x - k}{\epsilon\rho_1(n)})$ have 
disjoint support for integer 
$k$. Define
\begin{eqnarray*}
f_n(x) = \sum\limits_{k = -(n-1)}^{n-1}\phi(\frac {x - k}{\epsilon\rho_1(n)})
\end{eqnarray*}
and $$E_n = \bigcup\limits_{k = -(n-1)}^{n-1}[k - \epsilon\rho_1(n), k + \epsilon\rho_1(n)].$$ It 
is clear that $E_n$ is $\epsilon$-thin with respect to 
$\rho_1$, supp $f_n \subset E_n$ and $\|f_n\|_2^2 \sim n \epsilon\rho_1(n)$. We have
\begin{eqnarray*}
\hat f_n(y) = \sum\limits_{k = -(n-1)}^{n-1}e^{-2\pi i 
ky}\epsilon\rho_1(n)\hat{\phi}(\epsilon\rho_1(n)y) = D_{n-1}(y) 
\epsilon\rho_1(n)\hat{\phi}(\epsilon\rho_1(n)y).
\end{eqnarray*}
Pick an integer $a_n > 0$ which we will specify later. Let $$\s_n = \bigcup\limits_{l = -(a_n 
-1)}^{a_n -1}[l - \epsilon\rho_2(a_n), l + 
\epsilon\rho_2(a_n)]$$ 
then $\s_n$ is $\epsilon$-thin with respect to 
$\rho_2$. We have
\begin{eqnarray*}
\int\limits_{[-\frac 1 2, \frac 1 2] \backslash [ - \epsilon\rho_2(a_n),  \epsilon\rho_2(a_n)]} 
|D_{n-1}|^2 \lesssim \frac 1 { \epsilon \rho_2(a_n)}.
\end{eqnarray*}
Then 
\begin{eqnarray*}
\int\limits_{\s_n^c} |\hat f_n|^2 &\le& \int\limits_{\s_n^c \cap [-a_n, a_n]} |\hat f_n|^2 + 
\int\limits_{[-a_n, a_n]^c} |\hat f_n|^2 \\
 &\le& \sum\limits_{l = -(a_n -1)}^{a_n -1}\frac 1 { \epsilon\rho_2(a_n)}\epsilon^2\rho_1^2(n) + 
\sum\limits_{|l| \ge a_n} n \frac 
{\epsilon^2\rho_1^2(n)} {1 + (|l| \epsilon\rho_1(n))^{100}} \\
 &\lesssim& \frac {a_n \rho_1^2(n)\epsilon}{\rho_2(a_n)} + \frac {n\epsilon\rho_1(n)} {(a_n 
\epsilon\rho_1(n))^{99}}.
\end{eqnarray*}
Our goal is to make this expression much smaller than $\|f_n\|_2^2 \sim n \epsilon\rho_1(n)$. It 
will suffice to require that $1 \ll a_n \epsilon\rho_1(n)$ 
and $\frac {a_n \rho_1(n)}{n\rho_2(a_n)} \ll 1$. Set $C_1 = \frac k {\epsilon}$ and $C_2 = \frac 
{k^2} {\epsilon}$. Then there exists $t_k$ such that $\frac 
{C_2} 
{\rho_2(\frac 
{C_1} {\rho_1(t_k)})} < t_k$. By the way $t_k \rightarrow \infty$ as $k \rightarrow \infty$ since 
otherwise $\rho_2$ would be unbounded. 
Let $n = [t_k] > 0$. We also have that $\frac {C_1} {\rho_1(n)} \rightarrow \infty$ as $k 
\rightarrow \infty$ since otherwise $\rho_1$ would be 
unbounded. Let $a_n = \left [ \frac {C_1} {\rho_1(n)} \right ] > 0$. Note that $\rho_2(a_n) \ge 
\rho_2(\frac {C_1} {\rho_1(t_k)}) > \frac {C_2}{t_k}$. Then  
$a_n \epsilon\rho_1(n) \sim k$ and $\frac {a_n \rho_1(n)}{n\rho_2(a_n)} \lesssim \frac {k t_k} 
{\epsilon n C_2} \sim \frac 1 k$. Let $k \rightarrow \infty$ to 
obtain the desired result.\\

Now we will consider the case when dimension $d \ge 2$. The construction here is much simpler and 
we will just give a sketch. Let $\phi$ be a Schwartz function on $\Bbb{R}$ supported in $[-\frac 12, \frac 12]$. Define for $n > d$
$$f_n = 
\phi(\frac {x_1}{n 
- d})\prod\limits_{i = 2}^{d}\phi(\frac {x_i}{\epsilon\rho_1(n)})$$ and $$E_n = [-(n-d), 
n-d]\times[-\epsilon\rho_1(n),\epsilon\rho_1(n)]^{d-1}.$$ Then $E_n$ is $\epsilon$-thin with 
respect to $\rho_1$ and $f_n$ is supported in $E_n$. We have
$$\hat f_n (y) = (n - d)\hat \phi ((n-d)y_1)\prod\limits_{i = 2}^{d}\epsilon\rho_1(n)\hat 
\phi(\epsilon\rho_1(n)y_i).$$
Let $a_n > d$ be a number which we will specify later. Define 
$$\s_n = [-\epsilon\rho_2(a_n), \epsilon\rho_2(a_n)]\times[-(a_n - d), a_n - d]^{d - 1}$$ then 
$\s_n$ is $\epsilon$-thin with respect to 
$\rho_2$. We have
\begin{eqnarray*}
\int\limits_{\s_n^c} |\hat f_n|^2 \lesssim \frac {n\epsilon^{d-1}\rho_1^{d-1}(n)} {1 + 
(\epsilon\rho_2(a_n)n)^{100}} + \frac {n\epsilon^{d-1}\rho_1^{d-1}(n)}{1 
+ (a_n\epsilon\rho_1(n))^{100}}.
\end{eqnarray*}
This expression should be much smaller than $\|f_n\|_2^2 \sim n\epsilon^{d-1}\rho_1^{d-1}(n)$. 
Therefore it is enough to require that $\epsilon\rho_2(a_n)n \gg 
1$ and $a_n\epsilon\rho_1(n) \gg 1$. Set $C_1 = \frac k {\epsilon}$ and $C_2 = \frac {k} 
{\epsilon}$. Then there exists $n > d$ such that $\frac {C_2} 
{\rho_2(\frac 
{C_1} {\rho_1(n)})} < n$. Let $a_n = \frac 
{C_1} {\rho_1(n)} > d$ then $a_n\epsilon\rho_1(n) = k$  and $\epsilon\rho_2(a_n)n > \frac 
{\epsilon n C_2} {n} \ge k$. Let $k \rightarrow \infty$ to obtain 
the desired result.
\hfill$\square$\\


\begin{thebibliography}{99}

\bibitem{AB} W.O. Amrein, A.M. Berthier, {\sl On support properties of $L^p$ functions and
their Fourier trsnforms}, Journal of Functional Analysis, vol. 24, 1977, 258-267.

\bibitem{B} M. Benedicks, {\sl On Fourier transforms of functions supported
on sets of finite Lebesgue measure}, Journal of Mathematical Analysis and Applications,
vol. 106, 1985, 180-183.

\bibitem{FS} G. Folland, A. Sitaram, {\sl The Uncertainty Principle}, The Journal of
Fourier Analysis and Applications, vol. 3, 1997, 207-238.

\bibitem{HJ} V. P. Havin, B. Joricke, {\sl The Uncertainty Principle in Harmonic
Analysis},
Springer-Verlag, Berlin Heidelberg, 1994.

\bibitem{K1} O. Kovrijkine, {\sl Some results related to the Logvinenko-Sereda
theorem},
Proceedings of AMS, vol. 129, Num. 10, 2001, 3037-3047.

\bibitem{K3} O. Kovrijkine, {\sl Some estimates of Fourier transforms},
Ph.D.
thesis, Caltech, 2000.

\bibitem{K4} O. Kovrijkine, {\sl A version of the Uncertainty Principle for functions with 
lacunary Fourier transform},
preprint.

\bibitem{LS} V.N. Logvinenko, Yu. F. Sereda, {\sl Equivalent norms in spaces of entire functions 
of exponential type}, Teor. Funktsii, Funktsional. Anal. i 
Prilozhen 19, 1973, 234-246.

\bibitem{N} F. L. Nazarov, {\sl Local estimates of exponential polynomials
and their application to inequalities of uncertainty principle type},
St. Petersburg Math. J. 5(1994), 663-717.

\bibitem{R} W. Rudin, {\sl Complex and real analysis}, McGrow-Hill, 1990.

\bibitem{S} E. Stein, {\sl Harmonic analysis: real variable methods, orthogonality and 
oscillatory integrals}, Princeton University Press, 1993.

\bibitem{W} C. Shubin, R. Vakialian, T. Wolff, {\sl Some harmonic analysis 
questions suggested by Anderson-Bernoilli models},  Geom. Funct. Anal. 8 (1998), 
no. 5, 932--964.

\bibitem{Z} A. Zygmund, {\sl Trigonometric series, v. I and II}, Cambridge
University Press,
New York, 1968.

\end{thebibliography}
\end{document}